\documentclass{article}
\usepackage{graphicx} 
\usepackage{amsmath}
\DeclareMathOperator{\Arg}{Arg}
\usepackage{amssymb}
\usepackage{amsthm}
\usepackage{setspace}
\usepackage[margin=1in]{geometry}
\title{Numerical Shadow Rigidity for $3\times 3$ Matrices}
\author{Dylan Johnson \\ \small University of Iowa}
\date{April 2026}

\begin{document}
\newtheorem{thm}{Theorem}[section]
\newtheorem{lem}[thm]{Lemma}
\newtheorem{cor}[thm]{Corollary}
\newtheorem{prop}[thm]{Proposition}
\newtheorem{defi}[thm]{Definition}

\maketitle
\begin{abstract}
The \textbf{numerical shadow} of a complex matrix is the measure obtained by pushing forward that Haar measure on the unit sphere in $\mathbb{C}^n$ by the numerical range map. By the choice of measure, unitary equivalence implies same measure, but the reverse is not true for $n\geq 3$. We name a matrix $\textbf{shadow rigid}$ it its unitary equivalence class and shadow equivalence class coincide. It is known that normal matrices are rigid, but the topic has not been explored beyond this in the literature. The main result of the paper is showing that in the $n=3$ case, a matrix is rigid if and only if it is non-generic. We also show that rigidity is preserved under an affine group action and present three normal-like properties of $3\times 3$ non-generic matrices. 
\end{abstract}
\section{Introduction}
\begin{defi} Denoting the $n\times n$ complex valued matrices as $M_n(\mathbb{C})$ and the usual inner product on $\mathbb{C}^n$ as $\langle \cdot,\cdot \rangle$, we define the \textbf{numerical range} of a matrix $A$ as 
\begin{equation*}
    W(A):=\{\langle Ax,x\rangle\;|||x||=1\}.
\end{equation*} \end{defi}
Numerical ranges have been studied for over a century and there are at least two textbooks, \cite{yuangao} \cite{gusandrao}, on the subject.
\begin{defi} Given a matrix $A$ we define its \textbf{numerical range map}: 
\begin{align*}
    \phi_A: S_n&\to \mathbb{C}\\
    x&\mapsto \langle Ax,x\rangle,
\end{align*}
where $S_n$ denotes the unit sphere in $\mathbb{C}^n$. \end{defi} In \cite{numshadow} and \cite{nummeasure}, the authors define the \textbf{numerical shadow}, denoted $\mu_A$, of a matrix as the pushforward of the Haar measure on $S_n$, to the complex plane: 
\begin{equation*}
    \mu_A(U):=h(\phi_A^{-1}(U)),
\end{equation*}
where $U\subset \mathbb{C}$ is a measurable set and $h$ denotes the Haar measure, (also named Fubini-Study or Euclidean measure), on $S_n$. \\\\
By our choice of measure, it is obvious that unitarily equivalent matrices will give the same shadow. Less obvious is that the support of $\mu_A$ is equal to $W(A)$ and that the measure is absolutely continuous with respect to Lebesgue measure. Hence, we can define a density function for these measures. \\\\
See \cite{numshadow} and \cite{nummeasure} for an in-depth presentation of the details, but we give a few results to give the reader a flavor of these measures. For normal matrices, the numerical shadow is a multivariate B-spline with knots at its eigenvalues. For a $2\times2$ non-normal matrix, it is well known that the numerical range is a (possibly generate) ellipse. Let $E$ be the filled ellipse formed by the numerical range of a $2\times 2$ non-normal matrix, and $a,b$ give the lengths of the respective semimajor and semiminor axes, then the shadow density is given by 
\begin{equation*}
    \frac{1}{2\pi a b\sqrt{1-r^2}}
\end{equation*}
at every point on the elliptical curve $rE$.\\\\
In \cite{numshadow}, the authors note that as 
\begin{equation*}
    \langle Ax,x\rangle=tr(x^*Ax)=tr(Axx^*),
\end{equation*}
the numerical shadow may be viewed as the pushforward of a measure on the pure quantum states under the linear map $X\mapsto tr(AX)$. In another work \cite{numrangesareshadows}, the authors present a similar finding, showing that numerical range maps can be seen as an embedding followed by an affine map to the complex plane. \\\\
In \cite{geom}, the authors relate the numerical shadow to the geometry of state space. In particular, the authors show that all numerical ranges of $n\times n$ matrices are similar, up to scaling and translation, of orthogonal projections of $CP^{n-1}$ and that all orthogonal projections of $CP^{n-1}$ are realized as the numerical range of an $n\times n$ matrix. 
\begin{defi} For a smooth map between manifolds, we say a point in the domain is a \textbf{critical point} if the Jacobian does not have full rank at the point. We call a point in the range of the map a \textbf{critical value} if it is the image of a critical point.\end{defi}
The paper \cite{difftop}, provides a comprehensive analysis on the geometric properties of the numerical range map. Crucially for us, they prove that the critical points of the numerical range map of a matrix $A\in M_n(\mathbb{C})$ are exactly the eigenvectors of the parameterized family of Hermitian matrices $\Re(Ae^{-i\theta})$. In \cite{nummeasure} it is proven that the singular values of the shadow of a matrix are contained in the critical values of the numerical range map. They further conjecture that the sets are the same, and prove them to be the same in the two and three dimensional cases. 
\begin{defi} The authors define a matrix $A\in M_n(\mathbb{C})$ to be \textbf{generic} if the Hermitian matrices $\Re(Ae^{-i\theta})$ have $n$ distinct eigenvalues for all $\theta$.\end{defi} The authors prove that this is if and only if the numerical range maps of the paramaterized family are \textbf{Morse} for all $\theta$. \\\\
We list a few important lemmas from the paper. 
\begin{lem} The critical values of the numerical range map cross if and only if its critical points cross.\end{lem}
Proof: See \cite{difftop}.
\begin{lem} Non-generic matrices are Lebesgue measure zero in $M_n(C)$. \end{lem}
Proof: See \cite{difftop}
\begin{defi} We say that a square matrix is \textbf{shadow rigid} if its numerical shadow equivalence class is the same as its unitary equivalence class.\end{defi}
All $2\times2$ matrices are rigid. In fact, the stronger statement that $2\times 2$ matrices have the same numerical range if and only if unitarily equivalent is well known. In \cite{numshadow}, the authors prove that normal matrices are rigid for any dimension, and they give examples of matrices that are not similar, but have the same numerical shadow. Later in the paper, we will show that there are non-normal rigid $3\times 3$, matrices, but that rigid matrices form a Lebesgue measure zero set in $M_3(\mathbb{C})$; a sharp contrast to the $n=2$ case where all matrices are rigid.\\\\
We prove a handy lemma about rigidity preservation under a group action. \\\\
Let $G$ be the group of invertible affine transformations on $\mathbb{R}^2$. We define a group action on $M_n(\mathbb{C})$. Let $f\in G$ be defined by 
\begin{equation*}
    f\left(\begin{bmatrix}
        x_1\\ x_2
    \end{bmatrix}\right)=\begin{bmatrix}
        a_1 & b_1 \\ a_2 & b_2
    \end{bmatrix}\begin{bmatrix}
        x_1\\ x_2
    \end{bmatrix}+\begin{bmatrix}
        c_1 \\ c_2
    \end{bmatrix}.
\end{equation*}
Now let $A=\Re(A)+i\Im(A)\in M_n(\mathbb{C})$ be given. We define
\begin{equation*}
    f\cdot A:=(a_1\Re(A)+b_1\Im(A))+i(a_2\Re(A)+b_2\Im(A))+(c_1+ic_2)I_n.
\end{equation*}
\begin{lem} Rigidity is preserved under this affine group action. \end{lem}
Proof: The proof relies on the following two properties of this group action: \\
i) $f\circ\mu_A=\mu_{f\cdot A}$ and \\
ii) $f\cdot U^*AU=U^*f\cdot AU$ for any  unitary $U$ and $A\in M_n(\mathbb{C})$. \\\\
Now assume that $A$ is rigid and $f\in F$ is given; we show that $f\cdot A$ is rigid. Suppose $\mu_B=\mu_{f\cdot A}$. Then we have
\begin{align*}
    \mu_B&=\mu_{f\cdot A} \iff\\
    \mu_B&=f\cdot \mu_A \iff\\
    f^{-1}\mu_B=\mu_{f^{-1}\cdot B}&=\mu_A \iff\\
    f^{-1}\cdot B&=U^*AU \textbf{ (as $A$ is rigid)}\iff\\
    B=f\cdot U^*AU&=U^*(f\cdot A)U\quad \square
\end{align*}
In \cite{numshadow}, the following important theorem is proven: \\
\begin{thm}\label{trace} Two matrices $A,B\in M_n(C)$ give the same numerical shadow if and only if 
\begin{equation*}
    tr((sA+tA^*)^k)=tr((sB+tB^*)^k)
\end{equation*}
(as polynomials in $s,t$) for all $k\leq n$, which is if and only if 
\begin{equation*}
    \det(I_n-sA-tA^*)=\det(I_n-sB-tB^*)
\end{equation*}
as polynomials.\end{thm} 
We see that \textbf{Theorem \ref{trace}} resembles a weaker version of 
\begin{thm}\textbf{Specht's Theorem} 
 Two matrices $A,B\in M_n(\mathbb{C})$ are unitarily equivalent if and only if $tr(w(A,A^*))=tr(w(B,B^*))$ for all words $w$. \end{thm}
In light of the projective nature of the numerical range map and Specht's theorem, we have two ways of framing rigidity. First, rigid matrices highlight unique geometric features of $CP^{n}$. The image of the measure given by rigid matrices are obtained by a unique projection, mod unitary equivalence. On the other hand, like normal matrices, a rigid matrices' unitary class can be determined by this finite set of trace equations, going only up to the dimension of the matrix. \\\\
The rest of the paper is laid out as follows. In \textbf{Section \ref{examples}} we give examples of rigid and non-rigid matrices, with one family of examples leading us to conjecture the main result of the paper. In \textbf{Section \ref{main result}}, we show that a $3\times 3$ matrix is rigid if and only if it is non-generic. While generic matrices have come up frequently in work on numerical range maps, there doesn't appear to be any work on the algebraic properties of non-generic matrices. In \textbf{Section \ref{properties}}, we explore a few normal-like properties that non-generic $3\times 3$ matrices have. 
\section{Examples for $n=3$}\label{examples}
By cyclicity of trace, in the $3\times 3$ case, trace equations in \textbf{Theorem \ref{trace}} are all of the form $tr(w(A,A^*))$, where the word $w$ is of length less than or equal to three. Also note that as $tr(A^*)=\overline{tr(A)}$, only five of the equations are independent.\\\\
We see that equality $tr(A^k)_{k=1}^3$ means that two matrices with the same shadow must have the same eigenvalues. From equality of $tr(A^*A)$, they must have the same Frobenius norm. As the shadow is unitarily invariant, we work with upper triangular matrices. For $A=[a_{ij}]$ we have 
\begin{equation}\label{tr(A^2A)}
    tr(A^2A^*)=a_{11}(|a_{11}|^2+|a_{12}|^2+|a_{13}|^2)+a_{22}(|a_{22}|^2+|a_{12}|^2+|a_{23}|^2)+a_{33}(|a_{33}|^2+|a_{13}|^2+|a_{23}|^2)+a_{12}\bar{a_{13}}a_{23}.
\end{equation}
To show that a matrix is not rigid, we need to show a another matrix that is not unitarily equivalent, and gives the same trace values. Showing rigidity is a bit more difficult, but the problem is greatly simplified by restricting to upper triangular matrices. \\\\
In \cite{nummeasure} the authors show that the non-similar matrices $\begin{bmatrix}
    0 & 0 & 0 \\ 0 & 0 & 2 \\ 0 & 0 & 0
\end{bmatrix}$ and $\begin{bmatrix}
    0 & \sqrt{2} & 0 \\ 0 & 0 & \sqrt{2} \\ 0 & 0 & 0
\end{bmatrix}$ have the same numerical shadow. In fact the family $\begin{bmatrix}
    0 & \sqrt{4-4t} & 0  \\ 0 & 0 & 2\sqrt{t} \\ 0 & 0 & 0 
\end{bmatrix}_{0\leq t\leq 1}$ all give the same shadow. In light of the results in \cite{geom}, these equivalent shadows highlight some symmetry of state space. See that these examples break a plausible conjecture that the direct sum of rigid matrices would also be rigid. \\\\
One final note on the two non-rigid examples. Consider the map 
\begin{equation*}
    \begin{bmatrix}
        x_1\\x_2\\x_3
    \end{bmatrix}\mapsto \begin{bmatrix}
        1/\sqrt{2}(-x_1e^{-i2\arg(x_2)}+\bar{x_3})\\x_2\\1/\sqrt{2}(\bar{x_1}e^{i2\arg(x_2)}+x_3)
    \end{bmatrix}
\end{equation*}
This is injective and surjective on the unit sphere, and if we pre-compose it with the numerical range map of the first matrix, we get the numerical range map of the second. It is not clear if equivalent shadows implies the existence of such a map for $3\times 3$'s, but we know that it is true for the $n=2$ case. \\\\
As noted above, all normal matrices are rigid. An easy example of a rigid non-normal matrix to confirm in light of \textbf{Equation \ref{tr(A^2A)}} is $\begin{bmatrix}
    0 & 1 & 1 \\ 0 & 0 & 1 \\ 0 & 0 & 0
\end{bmatrix}$. Matrices with one eigenvalue are rigid if and only if they are unitarily equivalent to a matrix of the form $\begin{bmatrix}
    \lambda & re^{i\theta} & r \\ 0 & \lambda & r \\ 0 & 0 & \lambda
\end{bmatrix}$, where $r\geq 0$. It's not clear from this family what algebraic properties we would expect from a rigid matrix. We do see that no gaps in eigenvalues forces no gaps in the magnitude of the off diagonals.   \\\\
Matrices of the form $\begin{bmatrix}
    0 & t & 0 \\ 0 & 0 & 0 \\ 0 & 0 & 1
\end{bmatrix}$
are rigid for $0\leq t\leq 2$, but are not rigid for $t>2$. This family of examples led to a correct conjecture of the nature of rigidity for $3\times 3$'s. \\\\
What is going on here? It is known that the numerical range of a direct sum of matrices is the convex hull of the numerical ranges of the summands. The numerical range of $\begin{bmatrix}
    0 & t \\ 0 & 0 
\end{bmatrix}$ is the circle centered at the origin with radius $\frac{t}{2}$. So, for a fixed $t\geq 0$, the numerical range of our matrix is the convex hull of this circle and the point $\{1\}$. The $t=2$ case is where the eigenvalue $1$, lies on the boundary of the numerical range, and beyond this, it is properly contained inside the numerical range. \\\\
From our other rigid example, we know that eigenvalues residing on the boundary of the numerical range is not a necessary condition. We see what is going on when we examine the critical points of the numerical range map. Checking the eigenvalues of $\Re(A_t e^{-i\theta})$, we find $\{\pm \frac{t}{2},\cos(\theta)\}$. For $t\leq 2$, there exists a $\theta$ where the matrix $\Re(A_t e^{-i\theta})$ has a repeat eigenvalue, hence is non-generic. The main result of this paper is that for $3\times 3$ matrices, rigidity and non-genericity are equivalent. \\\\
For the previous two examples of rigid families, excepting the $t=2$ case, the authors in \cite{numrange3} prove the stronger result that the unitary classes of those matrices can be completely determined by their numerical range. They also prove that normal matrices whose eigenvalues are not colinear and those two families are the only classes of matrices that can be completely determined by the numerical range. In addition to these three classes, we will give one more class of matrices that can be determined by its shadow, and show that these  and normal matrices with colinear spectrums are the only classes that can be determined by the shadow.  
\section{Rigid if and only if Non-Generic}\label{main result}
\begin{thm}\label{mainthm} A $3\times 3$ matrix is rigid if and only if it is non-generic. \end{thm}
\subsection{Rigidity as Unique Pre-image}
In this subsection, we embed our unitary classes into a subset of $R^4$, and examine the image of this embedding under a function; equivalent to the mapping $A\mapsto tr(A^2A^*)$. The points of injectivity for this function represent the rigid matrices. \\\\
As noted in the previous section, four of the five trace equations fix the eigenvalues and Frobenius norm of our matrix. Now let $\{\lambda_i\}_{i=1}^3$ and $k\geq 0$ be given. Our $\lambda_i$'s will be eigenvalues and $k$ will be the Frobenius norm squared, less the magnitude squared of our eigenvalues. \\\\
By a theorem of Schur, every matrix satisfying our eigenvalue and norm criteria is unitarily equivalent to an upper triangular matrix with $a_{ii}=\lambda_i$.\\\\
We define the set
\begin{equation*}\label{gamma}
    \Gamma:=\{A\in M_3(C)\;|\; A=\begin{bmatrix}
        \lambda_1 & a_{12} & a_{13} \\ 0 & \lambda_2 & a_{23}\\ 0 & 0 & \lambda_3
    \end{bmatrix}\;\text{ with }\;|a_{12}|^2+|a_{13}|^2+|a_{23}|^2=k\}
\end{equation*}
We also define the subset of $R^3$
\begin{equation*}\label{D_k}
    D_k:=\{(u,v,w) \;|\; u,v,w\geq 0\;\&\; u+v+w=k\}.
\end{equation*}
We embed $\Gamma$ into $D_k\times [0,2\pi)$ like so:
\begin{equation*}
    A=[a_{ij}]\mapsto (|a_{12}|^2,|a_{13}|^2,|a_{23}|^2,\Arg(a_{12}\bar{a_{13}}a_{23})) 
\end{equation*}
where we set $\Arg(0)=0$. 
 \\\\
The astute reader has noticed that while this map is well defined, there are some issues with well-definedness for unitary classes when our matrices have degenerate eigenspaces. We acknowledge that the same equivalence class may have multiple representatives in $D_k\times [0,2\pi)$, but that this only happens on the boundary of $D_k$, (when one of the  upper off diagonal entries of the matrix are zero), and is irrelevant to the following arguments.\\\\
Facts: Denoting $D^\circ$ as the interior of $D_k$, every point in $D_k^\circ\times [0,2\pi)$ corresponds to an upper triangular matrix under this embedding, each unitary class has an image in $D_k\times [0,2\pi)$, and each distinct point in $D_k^\circ\times [0,2\pi)$ represents a distinct unitary class. \\\\
Now consider the function 
\begin{align*}
    F:D_k\times [0,2\pi)&\to R^2\\
    (u,v,w,\theta)&\mapsto (\lambda_1+\lambda_2)u+(\lambda_1+\lambda_3)v+(\lambda_2+\lambda_3)w+\sqrt{uvw}\cdot e^{i\theta},
\end{align*}
Combined with our embedding, this is the map from $\Gamma$ to $R^2$, as $A\mapsto tr(A^2A^*)-\sum_{i=1}^3\lambda_i|\lambda_i|^2$.\\\\
Occasionally, we will informally define $F$ directly on a matrix in $\Gamma$ using the notation $F(A)$ to simplify the discussion. We mean apply $F$ to the image of $A$ under this embedding.
\begin{lem}\label{injective} The matrix $A$ is rigid if and only if and only if $F$ is injective at $A$. \end{lem}
Proof: This lemma is clear besides the ambiguity of the same unitary equivalence class potentially having multiple representatives in $D_k\times[0,2\pi)$, (because of multiple upper triangular representatives). One remedy would be to do the same embedding, but delete an open edge or two, (depending on whether there are two or one eigenvalues for the matrices in $\Gamma$), from $D_k$, and restrict the domain of $F$ to this. Thus, ensuring that unitary classes have exactly one representative in our domain. We didn't do this explicitly to keep the argument short, and we will later see that we never have one of this ambiguous cases be rigid, so this domain deletion is not necessary. $\square$\\\\
We define a class of functionals on $\mathbb{R}^2$: 
\begin{align*}
    p_\theta: \mathbb{R}^2&\to \mathbb{R}\\
    \vec{v}&\mapsto \left\langle\vec{v},\begin{bmatrix}
        \cos(\theta)\\\sin(\theta)
    \end{bmatrix}\right\rangle.
\end{align*}
\begin{prop}\label{rigidprop} A matrix is rigid if and only if it maximizes $p_\theta\circ F$ for some $\theta$. \end{prop}
To prove this, we will show that the range of $F$ is convex and that the points of injectivity lay on the boundary of the image of $F$. The proposition follows from a fundamental property in convex analysis. 
\subsection{The Range of F is Convex}
\begin{lem}\label{convex} The range of $F$ is convex.\end{lem}
We prove a few supporting lemmas to show lemma \ref{convex} is true. \\\\
We write $L(u,v,w)=(\lambda_1+\lambda_2)u+(\lambda_1+\lambda_3)v+(\lambda_2+\lambda_3)w$ for the linear part of $F$.
\begin{lem}\label{disc} For each $(u,v,w)$ in $D_k$, the closed disc centered at $L(u,v,w)$ of radius $\sqrt{uvw}$ is contained in the range of $F$. \end{lem}
Proof: If the point is on the boundary of the closed disc, then clearly it is in the range of $F$. If it is in the interior, we can fix the $u$ and continuously decrease $v$, (increasing $w$ to stay in $D_k$), until our point is on the boundary of the circle centered at $L(u,v-t,w+t)$ with radius $\sqrt{u(v-t)(w+t)}$. $\square$.
\begin{lem}\label{strict} The function $r(u,v,w):=\sqrt{uvw}$ is strictly concave on $D_k$. \end{lem}
Proof: To reduce computations, we make the substitution $w=k-u-v$ and work with the function $r(u,v):=\sqrt{(k-u-v)uv}$.
We will use these partials and  Hessian argument 
\begin{align*}
    \frac{\partial^2}{\partial u^2}r(u,v)&=-\frac{v^2(k-v)^2}{4((k-u-v)uv)^{3/2}}<0\\
    \frac{\partial^2}{\partial v^2}r(u,v)&=-\frac{u^2(k-u)^2}{4((k-u-v)uv)^{3/2}}<0\\
    \frac{\partial^2}{\partial u\partial v}r(u,v)&=-\frac{uv(-k^2+3k(u+v)-2u^2-3uv-2v^2}{4((k-u-v)uv)^{3/2}}
\end{align*}
Looking at the determinant of the Hessian:
\begin{align*}
    r_{uu}r_{vv}-r_{uv}^2=\frac{u^2v^2(k-(u+v))^2(k(u+v)-(u+v)^2+uv)}{((k-u-v)uv)^{3/2}}>0
\end{align*}
on the interior of $D_k$. As the function is strictly concave on the interior of a compact and convex domain, it is concave on the whole domain. $\square$ 
\\\\
Using the same $L(u,v,w)$ from above, we consider the function
\begin{align*}
    f:R^2\times D_k&\to R\\
    (x,u,v,w)&\mapsto ||x-L(u,v,w)||-r(u,v,w).
\end{align*}
Fact: The function $f$ is convex. \\
Proof: The linear function ($x-L(u,v,w)$) is convex. Norms are convex, and the composition of convex functions is again convex. By the concavity of $r(u,v,w)$, its negation is convex. Finally, the sum of convex functions is again convex. $\square$ 
\\\\
Next we define 
\begin{equation*}
    C\subset R^2\times D_k,\quad C:=\{(x,u,v,w)\;|\;f(x,u,v,w)\leq 0\}.
\end{equation*}
Now as $C$ is a sub-level set of a convex function, it is a convex set (see \cite{convex}). This is what we need to complete the proof of \textbf{Proposition \ref{convex}}. Considering \textbf{Lemma \ref{disc}}, we see that the range of $F(u,v,w,\theta)$ is the $R^2$ component of $C$. As this is a projection of a convex set, it is convex. $\square$ 
\subsection{Rigid if and only if Boundary}
\begin{lem}\label{boundary}A matrix is rigid if and only if $F(A)$ maps to the boundary of the range of $F$.\end{lem}
Proof: We first show that points on the interior of the range of $F$ have multiple pre-images. \\\\
Let $z$ be a point in the interior of the range of $F$. There exists $(u_0,v_0,w_0)\in D_k$ such that $z$ is contained in the open disc centered at $L(u_0,v_0,w_0)$ with radius $\sqrt{u_0v_0w_0}$. Within $D_k$, we fix $u_0$ and let $v_0$ continuously approach zero, eventually we have $||L(u_0,v_0-t,w_0+t)-z||=\sqrt{u_0(v_0-t)(w_0+t)}$, and hence there exists a $\theta$ such that $F(u_0,v_0-t,w_0+t,\theta)=z$. Repeating this procedure, but letting $w_0$ go to zero, we obtain a different point in $D_k\times [0,2\pi)$ that maps to our interior point $z$. \\\\
Next, we show that boundary points in the range have unique pre-images. Suppose that a point \\$F(u,v,w,\theta)=z$ is on the boundary of the range of $F$. As this is a convex set, (\textbf{Lemma \ref{convex}}), there exists an angle $\tau$ such that $z$ maximizes $\left\langle \begin{bmatrix} z_x\\z_y\end{bmatrix},\begin{bmatrix}
    \cos(\tau)\\ \sin(\tau)
\end{bmatrix}\right\rangle$ over the range of $F$. See Theorem 11.1 in \cite{convex}. \\\\
We see that 
\begin{equation*}
    \left\langle \begin{bmatrix} z_x\\z_y\end{bmatrix},\begin{bmatrix}
    \cos(\tau)\\ \sin(\tau)
\end{bmatrix}\right\rangle=\Re((\lambda_1+\lambda_2)e^{-i\tau})u+\Re((\lambda_1+\lambda_3)e^{-i\tau})v+\Re((\lambda_2+\lambda_3)e^{-i\tau})w+\sqrt{uvw}\cos(\theta-\tau).
\end{equation*}
This equation is bounded above by the case when $\tau=\theta$, and as $z$ maximizes this inner product, we must have $\tau=\theta$ or $\sqrt{uvw}=0$. Restating this, $p_\tau\circ F(u,v,w,\theta)$ is maximized only if $\tau=\theta$ or $\sqrt{uvw}=0$. \\
\\
Fact: $p_\theta \circ F$ is stictly concave on the interior of $D_k\times\{\theta\}$.\\
Proof: On $D_k\times \{\theta\}$ we have
\begin{equation*}
    p_\theta\circ F=p_\theta\circ(\text{Linear + Radial}\cdot e^{i\theta})=\left\langle L(u,v,w),\begin{bmatrix}
        \cos(\theta)\\ \sin(\theta)
    \end{bmatrix}\right\rangle+\sqrt{uvw}
\end{equation*}
which is linear plus the strictly concave $r(u,v,w)$. $\square$\\\\
With this, we have that $p_\theta\circ F$ is a strictly concave function from a compact convex domain to $R$, hence it has a unique maximizing domain point. This completes the proof of \textbf{Lemma \ref{boundary}} $\square$. \\\\
This also implies that every point on the boundary of the image of $F$ is an extreme point. Meaning, the only boundary points of $D_k$ that could be boundary points in the image are  those where one of the components equals $k$, ($u=k$ for instance). The issue we had with one equivalence class mapping to multiple points on the boundary is now resolved as well.\\\\
The preceding argument along with \textbf{Lemma \ref{injective}} prove \textbf{Proposition \ref{rigidprop}}. $\square$
\subsection{Equations for Non-Generic}
\begin{prop}\label{genericprop} An upper triangular matrix is non-generic if and only if $p_\theta \circ F$ is maximized for some $\theta$. \end{prop}
We use the notation
\begin{equation*}
    A_\theta:=\Re(Ae^{-i\theta})=\frac{1}{2}(Ae^{-i\theta}+A^*e^{i\theta}).
\end{equation*}
For $A\in\Gamma$, we have the following trace equations for $A_\theta$:
\begin{align*}
    tr(A_\theta)&=\sum \Re(\lambda_ie^{-i\theta})\\
    tr(A_\theta^2)&=(\sum \Re(\lambda_ie^{-i\theta})^2)+\frac{1}{2}\sum_{j>i}|a_{ij}|^2\\
    tr(A_\theta^3)&=(\sum\Re(\lambda_ie^{-i\theta})^3)\\
    &+\frac{3}{4}(|a_{12}|^2(\Re(\lambda_1e^{-i\theta})+\Re(\lambda_2e^{-i\theta}))+|a_{13}|^2(\Re(\lambda_1e^{-i\theta})+\Re(\lambda_3e^{-i\theta}))+|a_{23}|^2(\Re(\lambda_2e^{-i\theta})+\Re(\lambda_3e^{-i\theta})))\\
    &+\frac{3}{4}|a_{12}||a_{13}|a_{23}|\cos(\tau-\theta)
\end{align*}
where $\tau=\arg(a_{12}\bar{a_{13}}a_{23})$.\\\\
We pause to highlight that the above are the evaluations of the trace polynomials in \textbf{Theorem \ref{trace}}, evaluated as $s=e^{-i\theta}$, and $t=\bar{s}$. \\\\
Rearranging the $tr(A_\theta^3)$ equation we see
\begin{equation*}
   \frac{4}{3} (tr(A_\theta^3)-(\sum\Re(\lambda_ie^{-i\theta})^3))=p_\theta\circ F(A)
\end{equation*}
To prove \textbf{Proposition \ref{genericprop}}, we first assume that $A$ is non-generic, manifested at $\theta$. The three trace equations for $A_\theta$ assign a value to $p_\theta\circ F(A)$. Over the next two subsections, we show that this value is indeed the maximum of $p_\theta \circ F(A)$, ranging over $A\in \Gamma$. By the unique solution set of $\{tr(A_\theta^k)\}_{k=1}^3$ and showing equality to this maximum, we also prove the reverse implication. Namely, maximizing this functional implies non-generic.\\\\
Assume that $A$ is non-generic Then there exists two real numbers $r_1,r_2$ and a rank 2 orthogonal projection $P$, and $\theta$ such that 
\begin{equation*}
    A_\theta=r_1 P+r_2(I_3-P).
\end{equation*}
Our trace equations give 
\begin{align*}
tr(A_\theta)&=2r_1+r_2\\
tr(A_\theta^2)&=2r_1^2+r_2^2\\
    tr(A_\theta^3)&=2r_1^3+r_2^3
\end{align*}
To simplify notation, I denote $\alpha_i=\Re(\lambda_ie^{-i\theta})$ and make the substitution $k=\sum_{j>i}|a_{ij}|^2$.\\\\
Solving the first two equations, we get 
\begin{align*}
    r_1&=\frac{\sum \alpha_i}{3}\pm\frac{\sqrt{4\sum \alpha_i^2-4\sum_{j>i}\alpha_i\alpha_j+3k}}{6}\\
    r_2&=\frac{\sum \alpha_i}{3}\mp\frac{\sqrt{4\sum \alpha_i^2-4\sum_{j>i}\alpha_i\alpha_j+3k}}{3}
\end{align*}
See that if $A_\theta$ is as above, then $A_{\theta+\pi}=-A_\theta$, hence we can pick $r_1=\frac{\sum \alpha_i}{3}-\frac{\sqrt{4\sum \alpha_i^2-4\sum_{j>i}\alpha_i\alpha_j+3k}}{6}$. We make this choice and plug into the $tr(A_\theta^3)$ equation: 
\begin{equation}\label{nongenericeq}
    2\left(\frac{\sum \alpha_i}{3}-\frac{\sqrt{4\sum \alpha_i^2-4\sum_{j>i}\alpha_i\alpha_j+3k}}{6}\right)^3+\left(\frac{\sum \alpha_i}{3}+\frac{\sqrt{4\sum \alpha_i^2-4\sum_{j>i}\alpha_i\alpha_j+3k}}{3}\right)^3-\sum\alpha_i^3=p_\theta \circ F(A)
\end{equation}
Next, we will use the Lagranian method to find the maximum value of $p_\theta\circ F$ on $D_k\times\{\theta\}$, and show that this maximum value agree with the above. 
\subsection{Finding function maximums}
As the functional $p_\theta\circ F(u,v,w,\tau)$ is bounded by the case when $\tau=\theta$, we assume this to be the case for this maximum. \\\\
We use the Lagrangian method to find the maximum of $p_\theta\circ F$. Explicitly,  
\begin{equation*}
    p_\theta\circ F(u,v,w,\theta)=\Re((\lambda_1+\lambda_2)e^{-i\theta})u+\Re((\lambda_1+\lambda_3)e^{-i\theta})v+\Re((\lambda_2+\lambda_3)e^{-i\theta})w+\sqrt{uvw}.
\end{equation*}
We write $a_1=\Re((\lambda_1+\lambda_2)e^{-i\theta}),\;a_2=\Re((\lambda_1+\lambda_3)e^{-i\theta}),\text{ and  }a_3=\Re((\lambda_2+\lambda_3)e^{-i\theta})$ to neaten up notation. See that the $a_i$'s can be written in terms of the $\alpha_i$ from the previous section. Later we will rewrite both in terms of elementary symmetric polynomials. \\\\
We define the constraint function 
\begin{equation*}
    g(u,v,w)=u+v+w-k=0.
\end{equation*}
With the Lagrange multiplier $\lambda$, our Lagrangian function is 
\begin{equation*}
    \Lambda(u,v,w,\lambda)=a_1u+a_2v+a_3w+\sqrt{uvw}-\lambda(u+v+w-k).
\end{equation*}
Taking partials, we have the following
\begin{align*}
    \frac{\partial \Lambda}{\partial u}&=a_1+\frac{1}{2}\sqrt{\frac{vw}{u}}-\lambda=0\rightarrow \lambda-a_1=\frac{1}{2}\sqrt{\frac{vw}{u}}\\    
        \frac{\partial \Lambda}{\partial v}&=a_2+\frac{1}{2}\sqrt{\frac{uw}{v}}-\lambda=0\rightarrow \lambda-a_2=\frac{1}{2}\sqrt{\frac{uw}{v}}\\
        \frac{\partial \Lambda}{\partial w}&=a_3+\frac{1}{2}\sqrt{\frac{uv}{w}}-\lambda=0\rightarrow \lambda-a_3=\frac{1}{2}\sqrt{\frac{uv}{w}}        
\end{align*}
Now we solve for $u,v,w$ in terms of $\lambda$, finding
\begin{align*}
    u&=4(\lambda-a_2)(\lambda-a_3)\\
    v&=4(\lambda-a_1)(\lambda-a_3)\\
    w&=4(\lambda-a_1)(\lambda-a_2).
\end{align*}
Substituting these into the constraint, we have 
\begin{align*}
    12\lambda^2-8(a_1+a_2+a_3)\lambda+4(a_1a_2+a_2a_3+a_3a_1)-k=0
\end{align*}
Then
\begin{equation*}
\lambda=\frac{2(a_1+a_2+a_3)\pm\sqrt{4(a_1+a_2+a_3)^2-12(a_1a_2+a_2a_3+a_1a_3)+3k}}{6}
\end{equation*}
And our maximums and minimums are 
\begin{equation*}
    (p_\theta\circ F)_{max/min}=8\lambda^3-4(a_1+a_2+a_3)\lambda^2+4a_1a_2a_3
\end{equation*}
Before we move on, we point out that we now have an explicit description of rigid matrices in $\Gamma$, using the $u,v,w$ equations above. In \textbf{Section \ref{properties}}, we give a second description using a basis determined by the orthogonal projection $P$. 
\subsection{Verifying Equality}
In this section we verify that $(p_\theta\circ F)_{max}$ and the left hand side of \textbf{Equation \ref{nongenericeq}} are equal. 
\\\\
We rewrite both expressions using the elementary symmetric polynomials in $\alpha_1,\alpha_2,\alpha_3$: 
\begin{align*}
    e_1&=\alpha_1+\alpha_2+\alpha_3\\
    e_2&=\alpha_1\alpha_2+\alpha_2\alpha_3+\alpha_1\alpha_3\\
    e_3&=\alpha_1\alpha_2\alpha_3
\end{align*}
As we would hope, the radical terms from both expression are equivalent:
\begin{equation*}
    D:=4(\sum a_i)^2-12(a_1a_2+a_1a_3+a_2a_3)+3k=4\sum\alpha_i^2-4\sum_{j>i}\alpha_i\alpha_j+3k=4(e_1^2-3e_2)+3k
\end{equation*}
Now rewriting \textbf{Equation \ref{nongenericeq}} terms of esp, we get
\begin{equation*}
    p_\theta\circ F(A)=\frac{-8}{27}e_1^3+\frac{4}{3}e_1e_2+\frac{2k}{3}e_1-4e_3+\frac{1}{27}D^{3/2}
\end{equation*}
From our Lagrangian equation, we have 
\begin{equation*}
    F_{max}=8\lambda^3-8e_1\lambda^2+4(e_1e_2-e_3),
\end{equation*}
and
\begin{equation*}
    \lambda=\frac{4e_1\pm\sqrt{D}}{6}.
\end{equation*}
Plugging this in and simplifying, we obtain
\begin{equation*}
    (p_\theta\circ F)_{max}=\frac{-32}{27}e_1^3+\frac{2}{9}D+4(e_1e_2-e_3)\pm \frac{D^{3/2}}{27}.
\end{equation*}
Now as $\frac{D^{3/2}}{27}\geq 0$, and this is a maximum, we must have 
\begin{equation*}
    (p_\theta\circ F)_{max}=\frac{-32}{27}e_1^3+\frac{2}{9}D+4(e_1e_2-e_3)+ \frac{D^{3/2}}{27}.
\end{equation*}
Taking the difference: 
\begin{align*}
    (p_\theta\circ F)_{max}- p_\theta\circ F(A)&=-\frac{8}{9}e_1^3+\frac{2}{9}e_1D+\frac{8}{3}e_1e_2-\frac{2k}{3}e_1\\
    &=e_1(-\frac{8}{9}e_1^2+\frac{2}{9}D+\frac{8}{3}e_2-\frac{2k}{3})\\
    &=e_1(-\frac{8}{9}e_1^2+\frac{2}{9}(4(e_1^2-3e_2)+3k)+\frac{8}{3}e_2-\frac{2k}{3})\\
    &=0.
\end{align*}
Now we did not assume that $\arg(a_{12}\bar{a_{13}}a_{23})=\theta$, but again as the functional is bounded by this case and we have shown it to be the maximum, we either have $\arg(a_{12}\bar{a_{13}}a_{23})=\theta$ or one of those terms is zero. \\\\
We have shown that if a matrix is non-generic, then it maximizes the functional $p_\theta\circ F$. For the other direction, if this functional is maximized, then the eigenvalues $\{r_1,r_1,r_2\}$, satisfy the trace equations $\{tr(A_\theta^k\}_{k=1}^3$, which have a unique solution. Hence, $A_\theta$ has a repeat eigenvalue, and $A$ is non-generic. This proves \textbf{Proposition \ref{genericprop}}. $\square$\\\\
Finally, \textbf{Proposition \ref{rigidprop}} and \textbf{Proposition \ref{genericprop}} give us our proof of \textbf{Thereom \ref{mainthm}}. $\square$
\section{Properties of Non-generic 3x3}\label{properties}
In this section we explore a few properties of of non-generic $3\times 3$ matrices. Restating the definition, if $A\in M_3(\mathbb{C})$ is non-generic, then there exists real numbers, $r_1,r_2,\theta$ and a rank-2 orthogonal projection $P$ such that 
\begin{equation*}
    Ae^{-i\theta}+A^*e^{i\theta}=r_1P+r_2(I_3-P).
\end{equation*}
To neaten the notation, we write $I_3-P=P^\perp$ for the remainder of this section.\\\\
One non-trivial observation is that $PAP$ is normal. To prove this, we borrow a few tools from \cite{difftop}. 
\begin{defi} A vector $v$ is called a \textbf{normal vector} for a matrix $A$ if it is an eigenvector of both $A$ and $A^*$. \end{defi}
\begin{lem}\label{normalvector} An eigenvector $v$ of $A\in M_n(\mathbb{C})$ with eigenvalue $\lambda$ is a normal eigenvector of $A$ if and only if $A$ is unitarily equivalent to $\begin{bmatrix}
    \lambda & 0 \\ 0 & A_2
\end{bmatrix}$, where $A_2\in M_{n-1}(\mathbb{C})$.\end{lem}
Proof: See \cite{difftop}
\begin{cor}\label{corr} A matrix is normal if and only if all of its eigenvectors are normal eigenvectors. \end{cor}
\begin{lem} If $A$ is as above, then $PAP$ is normal. Hence, $A$ has a rank-2 normal compression. \end{lem}
Proof: Taking \textbf{Equation \ref{nongenericeq}} and multiplying on both sides by $P$, we have 
\begin{equation*}
    PAPe^{-i\theta}+PA^*Pe^{i\theta}=r_1P.
\end{equation*}
Now as a non-zero vector in $ran(P^\perp)$ is in the kernel of both $PAP$ and $(PA^*P)^*=PA^*P$, it is a normal eigenvector.  Now using the decomposition given by \textbf{Lemma \ref{normalvector}}, $PAP$ has another eigenvector in the range of $P$. It is easy to check using the above equation, that this is also an eigenvector of $PA^*P=(PAP)^*$, hence it is a normal eigenvector. Applying the lemma again forces another normal eigenvector. Then by the \textbf{Corollary \ref{corr}}, $PAP$ is normal. $\square$\\\\
That is the first of three normal-like properties we have for non-generic $3\times 3$'s: they have a normal rank-2 compression. While this condition is necessary, it is not sufficient to be non-generic. \\\\
Going back to \textbf{Equation \ref{nongenericeq}}, then multiplying on the left by $P$ and the right by $P^\perp$, we have 
\begin{equation}\label{normaleq}
    PAe^{-i\theta}P^\perp+PA^*e^{i\theta}P^\perp=PAe^{-i\theta}P^\perp+(P ^\perp Ae^{-i\theta}P)^*=0.
\end{equation}
This gives our second normal-like property: with respect to the same $P$ where $PAP$ is normal, we have $||PAP^\perp||=||QAP^\perp||$.\\\\
Piecing these two lemmas together, we find that $A$ is unitarily equivalent to 
\begin{equation*}
    \begin{bmatrix}
        a_1 & 0 & b_1 \\ 0 & a_2 & b_2 \\ -\overline{b_1e^{2i\theta}} & -\overline{b_2e^{2i\theta}} & a_3
    \end{bmatrix}
\end{equation*}
for some complex numbers $\{a_i\},\{b_i\}$, and $\theta$ satisfying $\Re((a_1-a_2)e^{-i\theta})=0$. This is with respect to an orthonormal basis defined by $P$ and $P^\perp$
\begin{lem} If $A$ is as above, the commutator $[A^*,A]$ is rank-deficient. \end{lem}
Proof: It is sufficient to show that $P[A^*,A]P=0$ as $P^\perp [A^*,A]P$ has kernel. \\\\
See that from \textbf{Equation \ref{normaleq}}, we have that 
\begin{equation*}
    PAP^\perp=-PA^*P^\perp e^{2i\theta}.
\end{equation*}
A similar argument gives
\begin{equation*}
    P^\perp A^*P=-P^\perp A^*Pe^{-2i\theta}.
\end{equation*}
Putting these together, we have
\begin{equation*}
    PAP^\perp\cdot P^\perp A^*P=(-PA^*P^\perp e^{2i\theta})\cdot (-P^\perp APe^{-2i\theta})=PA^*P^\perp\cdot P^\perp AP.
\end{equation*}
Combining this with the normality of $PAP$ gives us the lemma: 
\begin{align*}
    P(A^*A-AA^*)P&=P(A^*(P+P^\perp)A-A(P+P^\perp)A^*)P\\
    &=PA^*PPAP-PAPPA^*P+PA^*P^\perp P^\perp AP-PAP^\perp P^\perp A^*P\\
    &=0.\quad \square
\end{align*}
\section{Acknowledgments}
The author gratefully acknowledges support from the Erwin and Peggy Kleinfeld Graduate Fellowship fund, support from a GAANN fellowship (U.S. Department of Education), and Dr. Palle. E.T Jorgensen for his guidance.

\end{document}